\documentclass{amsart}

\usepackage{amsthm}
\usepackage{enumerate}
\usepackage{amsmath,amssymb,amscd,amsfonts}
\usepackage{mathtools}
\allowdisplaybreaks[4]

\newtheorem{theorem}{Theorem}[section]
\newtheorem{lemma}[theorem]{Lemma}

\newtheorem{corollary}[theorem]{Corollary}

\theoremstyle{definition}

\theoremstyle{remark}
\newtheorem{remark}[theorem]{Remark}

\numberwithin{equation}{section}

\newcommand{\NN}{\mathbb{N}}

\newcommand{\RR}{\mathbb{R}}
\newcommand{\CC}{\mathbb{C}}
\newcommand{\DD}{\mathbb{D}}

\newcommand{\GG}{\mathcal{G}}
\newcommand{\FF}{\mathcal{F}}

\newcommand{\be}{\beta}

\newcommand{\ep}{\varepsilon}
\newcommand{\clos}{\operatorname{clos}}
\begin{document}

\title{Galois conjugates for some family of generalized beta-maps}



\author{Shintaro Suzuki}
\address{Department of Mathematics, 
Tokyo Gakugei University, 
4-1-1 Nukuikita-machi 
Koganei-shi, 
Tokyo 184-8501,
Japan}
\email{shin05@u-gakugei.ac.jp}

\subjclass[2020]{}

\maketitle

\begin{abstract}
A real number $\be>1$ is called an Yrrap (or Ito-Sadahiro) number if the corresponding negative $\be$-transformation defined by $x\mapsto 1-\{\be x\}$ for $x\in[0,1]$, where $\{y\}$ denotes the fraction part of $y\in\mathbb{R}$, has a finite orbit at $1$. 
Yrrap numbers are an analogy of Parry numbers for positive $\be$-transformations given by $x\mapsto \{\be x\}$ for $x\in[0,1]$, $\be>1$. In this paper, we determine the closure of the set of Galois conjugates of Yrrap numbers. In addition, we show an analogy of the result to the family of piecewise linear continuous maps each of which is obtained by changing the odd-numbered branches (left-most one is regarded as $0$-th) of the $\be$-transformation to negative ones for $\be>1$. 
As an application, we see that both the set of Yrrap numbers which are non-Parry numbers and that of Parry numbers which are non-Yrrap numbers are countable. 


\end{abstract}

\section{Introduction}
\label{intro}
Let  $\beta$ be a real number greater than $1$. The so-called $\be$-transformation $T_\be:[0,1]\to[0,1]$ is defined by
$$
T_\be(x)=\{\be x\}
$$
for $x\in[0,1]$, where $\{y\}$ denotes the fraction part of $y\in\RR$. 
In case  the map has a finite orbit at $1$, it is easily seen that $\be$ is an algebraic integer, which relates the dynamical properties of the map to the algebraic properties of $\be$. In that case, $\be$ is called a Parry number and such numbers have been studied from viewpoints of number theory and ergodic theory (\cite{It-Ta,Pa,Pa2,Re,So1}). 
One of the natural questions for Parry numbers is to ask which algebraic integers are actually Parry numbers. It is well-known that Pisot numbers are Parry numbers (see \cite{Sc}), although it has not been known if the same is true for all Salem numbers (see \cite{Bo1,Bo2}). 
As a direction to provide non-Parry numbers, Solomyak \cite{So1} determined the closure of the set of Galois conjugates of Parry numbers using zeros of power series whose coefficients are in $[0,1]$, which yields that real algebraic integers greater than $1$ whose Galois conjugates do not belong to that set are non-Parry numbers. In \cite {Th}, Thompson extended the result to the case of the family of all generalized $\be$-maps in the sense of G\'ora \cite{Go}, each of whose element is obtained by changing some branches of the $\be$-transformation by negative ones. 

In this paper, we investigate two subfamilies of generalized $\be$-maps in the sense of G\'ora. One is the family of negative $\be$-transformations defined by $x\mapsto 1-\{\be x\}$, $\be>1$, $x\in[0,1]$. The other is the family of piecewise linear continuous maps each of which is obtained by changing the odd-numbered branches (left-most one is regarded $0$-th) of the $\be$-transformation to negative ones for $\be>1$. 
The first main result in this paper is an analogy of the result by Solomyak to the family of  negative $\be$-transformations. 
We recall that $\be>1$ is called an Yrrap number if the corresponding negative $\be$-transformation has a finite orbit at $1$. In Theorem \ref{main A}, we determine the closure of the set of Galois conjugates of all Yrrap numbers, which gives the answer to an open question in \cite{Ma}. Together with the result by Solomyak \cite{So1}, our result yields that the closure of the set of Galois conjugates of Parry numbers and that of Yrrap numbers are symmetric with respect to the imaginary axis in the complex plane. As its application, we see that the modulus of Galois conjugates of Yrrap numbers is less than or equal to the golden ratio (Corollary \ref{golden}). In addition, we show that the set of Yrrap (resp. Parry) numbers which are non-Parry (resp. non-Yrrap) numbers is countable in Theorem \ref{main C}. 

The second main result is also an analogy of the result by Solomyak to the other family of piecewise linear continuous maps stated above (Theorem \ref{main B}). 
Comparing with the result for all generalized $\be$-maps by Thompson \cite{Th},
we see that the closure of the set of Galois conjugates corresponding to that family  coincides with that to the family of all generalized $\be$-maps, although the family of maps we consider is a special class of generalized $\beta$-maps.

This paper is organized as follows. In Section 2, we introduce basic notions  used throughout this paper and state the main results precisely. Section 3 is devoted to giving all the proofs of the main results.

\section{Preliminaries and Main results}
In the following, we denote by $\NN$ the set of all non-negative integers and 
assume that $\be>1$ is a non-integer. 

Let $E=(E(0), E(1), \dots,)\in\{0,1\}^\NN$ be an infinite 0-1 sequence. Set $I_{\be}(x)=[\be x]$ for $x\in [0,1]$, where $[y]$ denotes the integer part of $y\in\RR$. 
We define the piecewise linear map on the unit interval $\tau_{\be,E} :[0,1]\to[0,1]$ by 
$$\tau_{\beta,E}(x)=E(I_{\beta}(x))+(-1)^{E(I_{\beta}(x))}T_{\beta}(x)$$
for $x\in[0,1]$. Note that the definition of the map $\tau_{\be,E}$ is slightly different from 
the original one in that a 0-1 sequence in \cite{Go} is defined as a finite sequence with $([\be]+1)$-terms. In our definition, actually, the map $\tau_{\be,E}$ depends only on the first $([\be]+1)$-terms in $E$ for $\be>1$. The reason why we define a 0-1 sequence  as an infinite sequence is to consider a family of generalized $\be$-maps $\{\tau_{\be,E}\}_{\be>1}$, in which $E$ is defined independently of the integer part of $\be>1$. 

As in the case of $\beta$-transformations, we can consider an expansion of numbers in $[0,1]$ using the map $\tau_{\be,E}$ as follows. 
For $x\in[0,1]$ and every non-negative integer $n\geq0$, we define the digit function $d_{n}(\be, E, x)$ and the sign function $e_{n}(\be, E, x)$ by 
\[\begin{split}
&d_{n}(\beta,E,x)=E(I_{\beta}(\tau_{\be, E}^{n}(x)))+I_{\beta}(\tau_{\be, E}^{n}(x)), \\
&e_{n}(\beta,E,x)=(-1)^{E(I_{\beta}(\tau_{\be, E}^{n}(x)))},
\end{split}\]
and the cumulative sign function $s_{n}(\be, E, x)$ by
\[s_{n}(\beta,E,x)=
\begin{cases}
1, &n=0,\\
\prod_{i=0}^{n-1} e_{i}(\be, E, x), &n\geq 1.
\end{cases}\] 
By definition we can represent $x\in[0,1]$ by
\begin{align*}\label{expansion}
\begin{split}
x&=\frac{s_{0}(\be, E, x)d_{0}(\be, E, x)}{\beta}+\frac{s_{1}(\be, E, x)\tau_{\be, E}(x)}{\beta} \\
&=\frac{s_{0}(\be, E, x)d_{0}(\be, E, x)}{\beta}+\frac{s_{1}(\be, E, x)d_{1}(\be, E, x)}{\beta^2}+
\frac{s_{2}(\be, E, x)\tau_{\be, E}^2(x)}{\beta^2} \\
&=\cdots \\
&=\sum_{i=0}^{n-1}\frac{s_{i}(\be, E, x)d_{i}(\be, E, x)}{\beta^{i+1}}+\frac{s_{n}(\be, E, x)\tau_{\be, E}^{n}(x)}{\beta^{n}}
\end{split}
\end{align*}
for every positive integer $n\geq1$ (see Proposition 1 in \cite{Go}).
Since $\tau_{\be, E}^n(x)\in[0,1]$ and $s_{n}(\be, E, x)\in\{-1,1\}$ for $n\geq 0$, 
taking $n\to\infty$ in the right side of the above equation yields
\[x=\sum^{\infty}_{n=0}\frac{s_{n}(\be, E, x)d_{n}(\be, E, x)}{\beta^{n+1}}.\] 
We call the above series the $\tau_{\be,E}$-expansion of $x$ and the sequence of integers
$\{ s_{n}(\be, E, x)d_{n}(\be, E, x)\}_{n=0}^{\infty}$
the coefficient sequence of the $\tau_{\be, E}$-expansion of $x$.
Set  the monotone pieces of the map $\tau_{\be, E}$ as
\[J_{i}=
\begin{cases}
[\frac{i}{\beta},\frac{i+1}{\beta}) & (0\leq i \leq [\beta]-1)\\
[\frac{[\beta]}{\beta},1]           &  (i=[\beta])
\end{cases}
\]
and denote by $EP$ the set of all end points of $J_i$'s except $0$ and $1$, that is, set
\[EP=\Bigl\{\frac{1}{\beta}, \dots, \frac{[\beta]}{\beta} \Bigr\}.\]
If there exists a positive integer $n\geq1$ such that $\tau_{\be, E}^n(x)\in EP$, we have
\[x=\sum_{i=0}^{n-1}\frac{s_{i}(\be, E, x)d_{i}(x)}{\beta^{i+1}}+\frac{s_{n}(\be, E, x)n_0}{\beta^{n+1}},\]
where $n_0=\beta\tau_{\be, E}^n(x)\in\{1,\cdots, [\beta]\}$. In this case, we call the $\tau_{\be,E}$-expansion of $x$ finite. 

In the case of $\be$-transformations (resp.\ negative $\beta$-transformations), recall that $\be$ is called a Parry (resp.\ Yrrap) number if the set of the orbit of $1$ by the corresponding map is a finite set. In particular, $\be$ is called simple if the orbit of $1$ by the map eventually falls to $EP$ (see \cite{Th}, \cite{Li-St}).  Similarly, for $E\in\{0,1\}^\NN$, we call $\be$ a generalized Parry number with respect to (w.r.t) $E$ if the set of the orbit of $1$ by the map $\tau_{\be,E}$ is a finite set. In particular, we call $\be$ is simple w.r.t $E$ if the orbit of $1$ by the corresponding map eventually falls to $EP$.

Let $I\subset\RR$ be a bounded interval and let $\DD$ be the open unit disk in $\CC$. Set
\[\FF_I=\Bigl\{1+\sum_{n=1}^\infty a_n z^n\ ;\ a_i\in I \text{ for } i\geq1\Bigr\} \]
and 
\[\GG_I=\{z\in\DD\ ;\ \exists f\in\FF_{I} \text{ such that } f(z)=0\}.\]
For an infinite 0-1 sequence $E\in\{0,1\}^{\NN}$, denote by $\Phi(E)$ the closure of the set of Galois conjugates of all simple generalized Parry numbers w.r.t $E$ and by $\hat{\Phi}(E)$ the closure of the set of Galois conjugates of all generalized Parry numbers w.r.t $E$. 

We define the formal power series $\phi_{\be, E}(z)$ by 
\[\begin{aligned}
\phi_{\tau_{\be,E}}(z)=
\sum_{n=0}^{N-1}s_n(\be,E,1)d_n(\be,E,1)z^{n+1} 
                      +k_{0}s_{N}(\be,E,1)z^{N+1} 
\end{aligned}\]
if $\beta$ is simple, where $N=\min\{n\geq1 ;\ \tau_{\be,E}^n(1)\in EP\}$ and $k_{0}=\beta\tau^{N}_{\beta,E}(1)$, and by
\[\phi_{\tau_{\be,E}}(z)=\sum_{n=0}^{\infty}s_n(\be,E,1)d_n(\be,E,1)z^{n+1}\]
if $\beta$ is non-simple. Note that the convergence radius of this power series is greater than or equal to $1$ since the coefficient sequences 
$\{s_n(\be,E,1)d_n(\be,E,1)\}_{n=0}^\infty$ is bounded. 
In case that $\be$ is a generalized Parry number, the power series $\phi_{\tau_{\be,E}}$ can be extended to a rational function, and otherwise has the unit circle as its natural boundary by Szeg\"o's theorem, which states that if the coefficients $\{b_n\}_{n=0}^\infty$ of the power series $\sum_{n=0}^\infty b_nz^n$ line in some finite set in $\CC$, then either it has the unit circle as its natural boundary or $b_n$ is eventually periodic (see \cite{Sz}).

Set 
\[\Phi^*(E)=\clos(\{z^{-1}\in\CC\ ;\  \exists \be>1\text{ such that }1-\phi_{\tau_{\be,E}}(z)=0\}),\]
where clos$(A)$ denotes the closure of $A\subset\CC$. 
Since $1-\phi_{\tau_{\be,E}}(\be^{-1})=0$ for $\be>1$, it is easily seen that $\Phi(E)\subset\hat{\Phi}(E)\subset\Phi^*(E)$.

In the following, we set 
\[E_0=(0,0,0,\cdots),\ \  E_1=(1,1,1,\cdots) \text{ and } E_{alt}=(0,1,0,1,0,1,0, \cdots).\] Notice that $\{\tau_{\be,E_0}\}_{\be>1}$ is the family of all $\be$-transformations and $\{\tau_{\be,E_1}\}_{\be>1}$ is that of all negative $\be$-transformations. For $\be>1$ the map $\tau_{\be,E_{alt}}$ is the continuous piecewise linear map obtained by changing the odd-numbered branches (left-most one is regarded as $0$-th) of the $\be$-transformation by negative ones. For a set of complex numbers $A\subset\CC$, let $A^{-1}=\{z^{-1}\in\CC\ ;\ z\in A\}$ and $-A=\{-z\in\CC\ ;\ z\in A\}$. In \cite{So1}, Solomyak determine the closure of the set of all Galois conjugates of (simple) Parry numbers:
\begin{theorem}[Theorem 2.1 in \cite{So1}]\label{Soso}
\[\Phi(E_0)=\hat{\Phi}(E_0)=\Phi^*(E_0)=\clos(\DD)\cup\GG_{[0,1]}^{-1}.\]
\end{theorem}
\noindent
As its generalization, Thompson \cite{Th} determined the closure of the set of Galois conjugates of all generalized Parry numbers by extending Parry's criteria (Theorem 3 in \cite{Pa}) to the case of generalized $\be$-maps under some conditions.
\begin{theorem}[Theorem 5.7 in \cite{Th}]\label{Thth}
\[\bigcup_{E\in\{0,1\}^{\NN}}\Phi(E)=\bigcup_{E\in\{0,1\}^{\NN}}\hat{\Phi}(E)=\bigcup_{E\in\{0,1\}^{\NN}}\Phi^*(E)=\clos(\DD)\cup\GG_{[-1,1]}^{-1}.\]
\end{theorem}

The first main theorem in this paper is an analogy of the result by Solomyak to negative $\beta$-transformations:

\begin{theorem}\label{main A}
\[\Phi(E_1)=\hat{\Phi}(E_1)=\Phi^*(E_1)=\clos(\DD)\cup -\GG_{[0,1]}^{-1}.\]
\end{theorem}
\noindent
Together with Theorem \ref{Soso}, we can see that $\Phi(E_1)=-\Phi(E_0)$, which means that these sets are symmetric with respect to the imaginary axis in the complex plane.  

In \cite{So1} and \cite{Fl}, it was shown that $\GG_{[0,1]}$ includes the open disk whose radius is $(\sqrt{5}-1)/2$ and the value of its radius is best possible. This immediately implies  the following corollary: 
\begin{corollary}\label{golden}
If $z\in\CC$ is a Galois conjugate of an Yrrap number, then its modulus is less than or equal to the golden ratio $(1+\sqrt{5})/2$. 
\end{corollary}

The second main theorem in this paper is also an analogy of the result by Solomyak to the maps  
$\{\tau_{\be,E_{alt}}\}_{\be>1}$.

\begin{theorem}\label{main B}
\[\Phi(E_{alt})=\hat{\Phi}(E_{alt})=\Phi^*(E_{alt})=\clos(\DD)\cup\GG_{[-1,1]}^{-1}.\]
\end{theorem}
\noindent
Together with Theorem \ref{Thth}, this theorem shows that $\Phi(E_{alt})=\cup_{E\in\{0,1\}^\NN}\Phi(E)$, although $E_{alt}$ is just an element in $\{0,1\}^{\NN}$.


In \cite{Go}, G\'ora defined a generalization of Chebyshev maps
$F_{\beta} : [-1,1]\to[-1,1]$ by 
\[F_{\beta}(x)=\cos(\beta \arccos x)\]
for $x\in[-1,1]$, where $\beta>1$. If $\beta$ is a positive integer, the map is the well-known Chebyshev polynomial of $\be$-th order. He showed that it is topologically conjugate to the map $\tau_{\be,E_{alt}}$ (see Proposition 17 in \cite{Go}) and the $T_{\beta}$-invariant density is expressed as a function associated with the orbit $\{F_{\beta}^{n}(-1)\}_{n=0}^{\infty}$. As an application of Theorem \ref{main B}, we have the following result:

\begin{corollary}
The set of Galois conjugates of $\be$'s such that the orbit of $-1$ by $F_\beta$ is a finite set is equal to $\clos(\DD)\cup\GG_{[-1,1]}^{-1}$.
\end{corollary}

Finally, as a relevant work, we investigate the set of all Parry numbers and that of all Yrrap numbers. In \cite{Li-St}, Liao and Steiner gave two examples of $\be$'s such that one is a Parry number but a non-Yrrap number and the other is a Parry number but a non-Yrrap number (see Proposition 6.1 and 6.2 in \cite{Li-St}). As an application of Theorem \ref{Soso} and \ref{main A}, we construct countably many such 
$\be$'s. 

\begin{theorem}\label{main C}
Let $P$ be the set of all Parry numbers and let $Y$ be the set of all Yrrap numbers. Then 
$P\cap Y^c$ and $P^c\cap Y$ are countable sets. 
\end{theorem}


\section{Proofs of Main results}

This section is devoted to giving the proofs of the main results stated in the previous section. First, we provide the proof of Theorem \ref{main A} in the same analogy of the result for $\beta$-transformations by Solomyak via the following two lemmas. The first lemma, proved by Thompson in \cite{Th}, provides a sufficient condition for $\be$ to be a simple generalized Parry number. 

\begin{lemma}[Theorem 5.5 in \cite{Th}]\label{lemma 1}
For a positive integer $N\geq2$, suppose that 

\noindent
$M(0), \dots, M(N)$ are distinct non-zero integers such that 
$|M(i)|+1<M(0)$ for all $i\geq 1$ and $|M(j)|\neq|M(k)|-1$ for $j,k$ with $j\neq k$. Then the equation
\[1=\frac{M(0)}{x}+\cdots+\frac{M(N)}{x^{N+1}}\]
has a solution $\be>1$ which is a simple generalized Parry number.

\end{lemma}

\begin{remark}\label{note}
Note that the original theorem in \cite{Th} states that the solution $\be$ above is only a generalized Parry number. By Lemma 5.4 in \cite{Th}, this solution is in fact a simple generalized Parry number, i.e., $\tau_{\be,E}^N(1)\in EP$. Furthermore, by Lemma 5.3 in \cite{Th}, we know that $|M(n)|=d_n(\be,E,1)$ and $\text{sign}(M(n))=s_n(\be,E,1)$ for $0\leq n\leq N$, where $\text{sign}(x)$ denotes the function defined by $\text{sign}(x)=1$ if $x>0$ and $\text{sign} (x)=-1$ if $x<0$ for $x\in\RR$.
\end{remark}

The next lemma relates the generating function for the coefficient sequence of the $\tau_{\beta,E}$-expansion of $1$ to that for the orbit $\{\tau_{\be,E}^n(1)\}_{n=0}^\infty$. 

\begin{lemma}[Proposition 4.1 in \cite{Su}]\label{lemma 2}
For $z\in\CC$ with $|z|<1$, we have
\[1-\sum_{n=0}^N s_n(\be,E,1)d_n(\be,E,1)z^{n+1}=(1-\beta z)\sum_{n=0}^N s_n(\be,E,1)\tau_{\be,E}^n(1) z^n,\]
where $N$ is the minimal number such that $\tau_{\be,E}^N\in EP$, regarded as $+\infty$ if no such $N$ exists.
\end{lemma}

\noindent
{\bf Proof of Theorem \ref{main A}}\ \ 
We first see that $\Phi^*(E_1)\setminus\clos(\DD)\subset -\GG_{[0,1]}^{-1}$. By Lemma \ref{lemma 2}, we have that
\[1-\phi_{\be,E_1}(z)=(1-\be z)\Biggl(1+\sum_{n=1}^N (-1)^n\tau_{\be,E_1}^n(1)z^n\Biggr).\]
By the definition of 
$\Phi^*(E_1)$, if $z_0\in\Phi^*(E_1)\setminus\clos(\DD)$, then $z_0^{-1}$ is a convergent of zeros of $1-\phi_{\be,E_1}(z)$ for $\be>1$ except $1/\be$. This ensures that $z_0^{-1}$ is also a convergent of zeros of $1+\sum_{n=1}^N \tau_{\be,E_1}^n(1)(-z)^n$, which means that $-z_0^{-1}\in\GG_{[0,1]}$, together with $\tau_{\be,E_1}^i(1)\in[0,1]$ for $i\geq1$. That is, we have $z_0\in -\GG_{[0,1]}^{-1}$.

We next see that $-\GG_{[0,1]}^{-1}\subset\Phi(E_1)\setminus\clos(\DD)$. Take a polynomial $f(z)$ of the form $f(z)=1+\sum_{n=1}^N a_n z^n$, where $N\geq1$ and $a_n\in(0,1)$ for $1\leq n\leq N$ with $a_i\neq a_j$ if $i\neq j$. Note that it is sufficient to show that for every positive number $\ep>0$, there exists a simple $\be>1$ such that 

\begin{equation}\label{app}
|\tau_{\be, E_{1}}^{n}(1)-a_n|<\ep
\end{equation}
for $1\leq n \leq N$ and

\begin{align}\label{irr}
x^N-d_0(\be, E_1)x^{N-1}+&\cdots+(-1)^{N-1}d_N(\be,1) \\ 
&=z^N(1-\be/z)\Biggl(1+\sum_{n=1}^{N}\tau_{\be, E_1}(1)\Bigl(\frac{-1}{z}\Bigr)^n\Biggr) \notag
\end{align}
is an irreducible polynomial. In fact, the irreducibility of the polynomial (\ref{irr}) ensures that the reciprocal of every zero of $1+\sum_{n=1}^N\tau_{\be,E_1}^n(1)(-z)^n$ is a Galois conjugate of $\be$.
Furthermore, the inequality (\ref{app}) ensures that any coefficient $a_n$ of the power series $f(z)$ can be approximated by $\tau_{\be,E_1}^n(1)$ for simple $\be$, which implies that there is a sequence $\{\be_m\}_{m=1}^\infty$ of simple numbers such that $1+\sum_{n=1}^N\tau_{\be_m,E_1}(1)z^n$ converges to $f(z)$ uniformly on compact subsets of $\DD$ as $m\to\infty$. By Hurwitz's theorem, we have that every zero of $f(z)$ in $\DD$ is obtained as the limit of zeros of $1+\sum_{n=1}^N\tau_{\be_m,E_1}(1)z^n$
as $m\to\infty$. Since the set of all polynomials of the form $f(z)$ is dense in the set $\FF_{[0,1]}$ with uniform topology, we get the conclusion. 

Set $a_0=0$, $a_{N+1}=1$ and $\delta=\min_{0\leq i,j\leq N+1}|a_i-a_j|$. For a positive integer $M\geq1$, we can see that $I_\be(a_n)\in\{I_M(a_n), I_M(a_n)+1\}$ for $\be\in(M,M+1)$ and $1\leq n \leq N$. Take an odd integer $M$ so large as 
$M>\max\{7/\delta, 5/\ep\}$. Let $D(0)=M+1$. 
For $1\leq n\leq N-1$, set $D(n)=I_M(a_n)$ if $I_M(a_n)$ is even and $D(n)=I_M(a_n)+1$ if $I_M(a_n)$ is odd. Note that
each $D(n)$ is  an even integer with $0\leq D(n)-I_M(a_n)\leq 1$. 
Define the even integer $D(N)$ by $D(N)=I_M(a_N)+k$, where 
\[k=
\begin{cases}
2, & (I_M(a_N)\equiv 0\mod 4) \\
1, & (I_M(a_N)\equiv 1\mod 4) \\
0, & (I_M(a_N)\equiv 2\mod 4) \\
-1. & (I_M(a_N)\equiv 3\mod 4)
\end{cases}
\]
By definition, $D(N)\equiv2$ (mod $4$) and $|D(N)-I_M(a_N)|\leq2$. Note that by $M>7/\delta$, we have
\[|D(i)-D(j)|\geq|I_M(a_i)-I_M(a_j)|-3\geq M|a_i-a_j|-4\geq3\]
for $0\leq i, j\leq N$ with $i\neq j$. By Lemma \ref{lemma 1} and \ref{note}, the equation 
\[1=\frac{D(0)}{x}-\frac{D(1)}{x^2}+\cdots+(-1)^N\frac{D(N)}{x^{N+1}}\]
has a unique positive solution $\be_0\in(M,M+1)$, which is a simple Yrrap number with
\[d_n(\be_0,E_1,1)=D(n)\]
for $0\leq n \leq N$. 

Since every $D(n)$ is even for $0\leq n\leq N$ and  $D(N)\equiv2$ (mod $4$), the polynomial (\ref{irr}) for $\be_0$ is irreducible by Eisenstein's criterion. By definition, 
\[\frac{I_{\be_0}(\tau_{\be_0,E_1}^n(1))}{\be_0}\leq \tau_{\be_0,E_1}^n(1)<\frac{I_{\be_0}(\tau_{\be_0,E_1}^n(1))+1}{\be_0}\]
and
\[\frac{I_{\be_0}(a_n)}{\be_0}\leq a_n<\frac{I_{\be_0}(a_n)+1}{\be_0}\]
for $1\leq n\leq N$. This shows 
\[|\tau_{\be_0,E_1}^n(1)-a_n|\leq\frac{|I_{\be_0}(\tau_{\be_0,E_1}^n(1))-I_\be(a_n)|+1}{\be_0}.\]
Since $D(n)=d_n(\be_0,E_1,1)=I_{\be_0}(\tau_{\be_0,E_1}^n(1))+1$, $|D(n)-I_M(a_n)|\leq2$ and $I_{\be_0}(a_n)\in\{I_M(a_n), I_{M}(a_n)+1\}$ for $0\leq n \leq N$, we have
\[|I_{\be_0}(\tau_{\be_0,E_1}^n(1))-I_\be(a_n)|\leq|D(n)-I_M(a_n)|+2\leq4,\]
which yields
\[|\tau_{\be_0,E_1}^n(1)-a_n|\leq\frac{5}{\be_0}<\frac{5}{M}<\ep,\]
as desired.

Finally, we see that $\clos(\DD)\subset\Phi(E_1)$. Note that each zero of a polynomial of the form $f_r(z)=1+z/r+\dots+z^N/r^N$, where $N\geq2$ and $r\in(1,+\infty)$ can be approximated by Galois conjugates of elements in $\Phi(E_1)$ as in the above argument. Since the set of all zeros of $f_r(z)$ is equal to the set $\{r e^{2\pi in/N}\}_{n=0}^{N-1}$, we can conclude that $\clos(\DD)\subset\Phi(E_1)$ by the fact that the set  
$\cup_{r>1}\cup_{N\geq1}\{re^{2\pi in/N}\}_{n=0}^{N-1}$ is dense in $\{z\in\CC\ ;\ |z|\geq1\}$, which finishes the proof. \qed

\begin{remark}
In the second part of the above proof, we take a simple Yrrap number satisfying (\ref{app}) and (\ref{irr}) so that its integer part is odd. Then the proof in fact yields the stronger result: Denoting by $\Phi_{odd}(E_1)$ the closure of the set of all Galois conjugates of simple Yrrap numbers each of whose integer part is odd, we have $\Phi_{odd}(E_1)=\clos(\DD)\cup-\GG_{[0,1]}^{-1}$.

\end{remark}

Next we show Theorem \ref{main B} with some modifications of the above proof.

\vspace{0.1cm}

\noindent
{\bf Proof of Theorem \ref{main B}}\ \ 
As in the first part of the proof of Theorem \ref{main A}, the equation by Lemma \ref{lemma 2} 
\[1-\phi_{\be,E_{alt}}(z)=(1-\be z)\Biggl(1+\sum_{n=1}^{N}s_{n}(\be,E_{alt},1)\tau_{\be,E_{alt}}^n(1)z^n\Biggr),\]
where $N$ is the minimal number with $\tau_{\be,E_{alt}}^N(1)\in EP$ regarded as $+\infty$ if no such $N$ exists, yields that $\Phi^*(E_{alt})\setminus\clos(\DD)\subset\GG_{[-1,1]}^{-1}$. 

We shall show that  $\GG_{[-1,1]}^{-1}\subset\Phi(E_{alt})\setminus\clos(\DD)$. 
Let $f(z)$ be a polynomial of the form 
\[f(z)=1+\sum_{n=1}^N c_n a_n z^n,\]
where $N\geq2$, $c_n\in\{-1,1\}$ and $a_n\in (0,1)$ for $1\leq n\leq N$ with $a_i\neq a_j$ for $1\leq i<j\leq N$. As in the same reason for the proof of Theorem \ref{main A}, it is sufficient to show that for every positive number $\ep>0$, there exists a simple $\be_0>1$ such that 

\begin{equation}\label{app*}
s_{n}(\be_0,E_{alt},1)=c_n, \ \ 
|\tau_{\be, E_{1}}^{n}(1)-a_n|<\ep
\end{equation}
for $1\leq n \leq N$ and 

\label{irr*}
\begin{align}
&x^N-s_{0}(\be_0,E_{alt},1)d_0(\be_0, E_{alt})x^{N-1}-\cdots-s_{N}(\be_0,E_{alt},1)d_N(\be_0, E_{alt}, 1) \\
&=x^N(1-\be/x)(1+\sum_{n=1}^{N}s_{n}(\be_0,E_{alt},1)\tau_{\be_0, E_{alt}}^n(1)(1/x)^n \notag
\end{align}
is an irreducible polynomial. Set $a_0=0$, $a_{N+1}=1$ and $\delta=\min_{0\leq i,j\leq N+1}|a_i-a_j|$. 
Take an even integer $M$ so large as $M>\max\{7/\delta, 5/\ep\}$. Define the sequence of even integers $\{D(n)\}_{n=0}^N$ by $D(0)=M$,  
\[D(n)=\begin{cases}
I_M(a_n) & \text{ if } I_M(a_n) \text{ is even, } \\
I_M(a_n)+1 & \text{ if } I_M(a_n) \text{ is odd, } 
\end{cases} 
\]
for $1\leq n \leq N-1$ and $D(N)$ as the even positive integer with $D(N)\equiv2 \text{ (mod $4$) }$ with 
 $|D(N)-I_M(a_N)|\leq2$. Then we can see in the same way of the proof of Theorem \ref{main A} that the sequence $\{c_n D(n)\}_{n=0}^N  $ satisfies the assumptions of Lemma \ref{lemma 1}. Hence there are a simple $\be_0>1$ and some $0$-$1$ sequence $E'\in\{0,1\}^{\NN}$ such that the inequality (\ref{app*}) holds for $1\leq n\leq N$ and the polynomial (\ref{irr*}) for $(\be,E')$ is irreducible with
$d_{n}(\be_0,E',1)=D(n)$ and $s_{n}(\be_0,E',1)=c_n$ for $1\leq n \leq N$. 

In fact, we can replace $E'$ by $E_{alt}$ from the following argument. Since every $D(n)$ is even and $e_{n}(\be_0,E',1)=s_{n}(\be_0,E',1)s_{n+1}(\be_0,E',1)=c_n c_{n+1}$ for $1\leq n\leq N$ by definition, we have
\[I_\be(\tau_{\be_0,E'}^n(1))=D(n) \text{ is an even integer if } e_{n}(\be_0,E',1)=+1\]
and
\[I_\be(\tau_{\be_0,E'}^n(1))=D(n)-1 \text{ is an odd integer if } e_{n}(\be_0,E',1)=-1\]
for $0\leq n\leq N$. This means that $E(\tau_{\be_0,E'}^n(1))=0$ if $I_{\be}(\tau_{\be_0,E'}^n(1))$ is even and 

\noindent
$E(\tau_{\be_0,E'}^n(1))=1$ if $I_{\be}(\tau_{\be_0,E'}^n(1))$
is odd. Therefore, we conclude that $\tau_{\be_0,E'}^n(1)=\tau_{\be_0,E_{alt}}^n(1)$ and
$s_{n}(\be_0,E',1)=s_n(\be_0,E_{alt},1)$ for $1\leq n\leq N$, which allows us to replace $E'$ by $E_{alt}$.

We can show $\clos(\DD)\subset\Phi(E_{alt})$ by the same argument as the last part of the proof of Theorem \ref{main A}. \qed

Finally, we prove Theorem \ref{main C} as an application of Theorem \ref{main A}.
\vspace{0.1cm}

\noindent
{\bf Proof of Theorem \ref{main C}}\ \ 
Let $f\in\FF_{[0,1]}$. By the definition of $\FF_{[0,1]}$, the power series $f$ has the form $f(z)=1+\sum_{n=1}^\infty a_n z^n$, where $a_n\in[0,1]$ for $n\geq1$. Since $f$ has the convergence radius greater than or equal to $1$ and $a_n\in[0,1]$ for $n\geq1$, we have that $f(t)=1+\sum_{n=1}^\infty a_n t^n\geq1>0$ for $0<t<1$, which means that $\GG_{[0,1]}\cap(0,+\infty)=\emptyset$. Together with Theorem \ref{Soso}, we have that every Parry number has no positive Galois conjugate greater than $1$. In addition, we note that the set $\GG_{[0,1]}$ has a `spike' on the negative real axis, which means that $\GG_{[0,1]}\cap(-1,0]=(-1, -(\sqrt{5}-1)/2]$ and there is $0<C<(\sqrt{5}-1)/2$ such that every $z\in\GG_{[0,1]}\setminus(-1,0]$ satisfies $|z|<C$ (see Corollary 2.3 and Lemma 4.1 in \cite{So1} ). Since $\Phi(E_1)\supset-\GG_{[0,1]}^{-1}$ by Theorem \ref{main A}, there are countably many distinct Yrrap numbers one of whose Galois conjugates lies in $(1, (\sqrt{5}+1)/2)$. Together with $\GG_{[0,1]}\cap(0,\infty)=\emptyset$, we have that they are non-Parry numbers, ensuring the set $P^c\cap Y$ is countable. Similarly, $\Phi(E_0)\supset\GG_{[0,1]}^{-1}$ by Theorem \ref{Soso}, there are countably many distinct Parry numbers one of whose Galois conjugates lies in $(-(\sqrt{5}+1)/2, -1)$. Together with $-\GG_{[0,1]}\cap(-\infty,0)=\emptyset$, we have that they are non-Yrrap numbers, ensuring the set $P\cap Y^c$ is countable.\qed

\vspace{1.5mm}

Acknowledgement:
This work was supported by JSPS KAKENHI Grant Number 20K14331.

\end{document}